\documentclass[11pt]{amsart}
\usepackage{amssymb}

\newtheorem{lemma}{Lemma}[section] 

\newtheorem{theorem}[lemma]{Theorem}

\newtheorem{proposition}[lemma]{Proposition}
\newtheorem{corollary}[lemma]{Corollary}

\begin{document} 
\title{Euler characteristic and Lipschitz-Killing curvatures of closed semi-algebraic sets}
\author{Nicolas Dutertre}
\address{Universit\'e de Provence, Centre de Math\'ematiques et Informatique,
39 rue Joliot-Curie,
13453 Marseille Cedex 13, France.}
\email{dutertre@cmi.univ-mrs.fr}

\thanks{
 Supported by {\em Agence Nationale de la Recherche} (reference ANR-08-JCJC-0118-01)}

\begin{abstract}
We prove a formula that relates the Euler-Poincar\'e characteristic of a closed semi-algebraic set to its Lipchitz-Killing curvatures.
\end{abstract}

\maketitle
\markboth{Nicolas Dutertre}{Euler characteristic and Lipschitz-Killing curvatures of closed semi-algebraic sets }

\section{Introduction}
In [Fu2], Fu developed integral geometry for compact subanalytic sets. Using methods from geometric measure theory, he associated with every compact subanalytic set $X \subset \mathbb{R}^n$ a sequence of curvature measures:
$$\Lambda_0(X,-), \ldots ,\Lambda_n(X,-),$$
called the Lipschitz-Killing measures and established several integral geometry formulas. Among them, he proved the following Gauss-Bonnet formula: $\Lambda_0(X,X)=\chi(X)$.
Another approach based on stratified Morse theory was presented later by Broecker and Kuppe [BK].

In [Du], we extended the Gauss-Bonnet formula to closed semi-algebraic sets. Let us describe it. Let $X \subset \mathbb{R}^n$ be a closed semi-algebraic set and let $(K_R)_{R>0}$ be an exhaustive family of compact subsets of $X$. The limit $\lim_{R \rightarrow +\infty} \Lambda_0(X,X\cap K_R)$ is finite and independent on the choice of the family. We denote it by $\Lambda_0(X,X)$. In [Du] Corollary 5.7, we proved:
$$\Lambda_0(X,X)= \chi(X)-\frac{1}{2} \chi(Lk^\infty(X)) 
\hfill -\frac{1}{2 \hbox{vol}(S^{n-1})} 
\int_{S^{n-1}} \chi(Lk^\infty(X \cap v^\perp))dv,$$
where $Lk^\infty (X)=X \cap S_R^{n-1}$, $R \gg 1$, is the link at infinity of $X$. 

In this paper, following a suggestion of Michel Coste, we continue the study of the Lipschitz-Killing curvatures of closed semi-algebraic sets. Before stating our main results, we need some notations:
\begin{itemize}
\item for $k \in \{0,\ldots,n \}$, $G_{n}^k$ is the Grassmann manifold of $k$-dimensional linear
subspaces in $\mathbb{R}^{n}$ and $g_n^k$ is its volume,
\item for $k \in \mathbb{N}$, $b_k$ is the volume of the $k$-dimensional unit ball and $s_k$ is the
volume of the $k$-dimensional unit sphere,
\item  for $R>0$, $B_R$ (resp. $S_R$) will denote the ball (resp. the sphere) centered at the origin of radius $R$ in $\mathbb{R}^n$.
If $R=1$, we will write $B^n$ (resp. $S^{n-1}$).
\end{itemize}
Let $X \subset \mathbb{R}^n$ be a closed semi-algebraic set. We will prove (Theorem 3.7 and Corollary 3.8):
$$\displaylines{
\qquad \lim_{R \rightarrow + \infty} \frac{\Lambda_k(X,X\cap B_R)}{b_k R^k} =- \frac{1}{2g_n^{n-k-1}} \int_{G_n^{n-k-1}} \chi(Lk^\infty(X \cap H)) dH \hfill \cr
\hfill + \frac{1}{2g_n^{n-k+1}} \int_{G_n^{n-k+1}} \chi(Lk^\infty(X \cap L)) dL\hbox{ for } k \in \{1,\ldots,n-2\}, \qquad  \cr
}$$
and:
$$\lim_{R \rightarrow + \infty} \frac{\Lambda_{n-1}(X,X\cap B_R)}{b_{n-1} R^{n-1}} = \frac{1}{2g_n^2} \int_{G_n^2} \chi(Lk^\infty(X \cap H)) dH,$$
$$\lim_{R \rightarrow + \infty} \frac{\Lambda_{n}(X,X\cap B_R)}{b_{n} R^{n}} = \frac{1}{2g_n^1} \int_{G_n^1} \chi(Lk^\infty(X \cap H)) dH.$$
One should mention that similar integral geometric formulas were obtained for closed subanalytic germs in [CM].

Combining this with the above Gauss-Bonnet formula, we will obtain (Theorem 3.9):
$$\chi(X)= \Lambda_0(X,X) + \sum_{k=1}^n    \lim_{R \rightarrow + \infty} \frac{\Lambda_k(X,X \cap B_R)}{b_k R^k}.$$
When $X$ is smooth of dimension $d$, $0<d<n$, this last formula can be reformulated in the following form (Theorem 4.3). If $d$ is even, then we have:
$$\chi(X)= \frac{1}{s_{n-1}} \int_X K_d(x)dx + \sum_{i=0}^{\frac{d-2}{2}} \lim_{R \rightarrow + \infty} \frac{1}{s_{n-d+2i-1}b_{d-2i}R^{d-2i}} \int_{X \cap B_R} K_{2i}dx.$$
If $d$ is odd, then we have:
$$\chi(X)=  \sum_{i=0}^{\frac{d-1}{2}} \lim_{R \rightarrow + \infty} \frac{1}{s_{n-d+2i-1}b_{d-2i}R^{d-2i}} \int_{X \cap B_R} K_{2i}dx.$$
Here the $K_i$'s are the classical Lipschitz-Killing-Weyl curvatures. 

The paper is organized as follows: in Section 2, we give the definition of the Lipschitz-Killing measures. In Section 3, we prove our main formulas, first in the conic case and then in the general case. In Section 4, we focus on the smooth case. We start recalling some Gauss-Bonnet formulas, due to several authors, for complete manifolds. Then, applying the results of Section 3, we generalize and improve these formulas in the case of smooth closed semi-algebraic sets.

\section{Lipschitz-Killing curvatures of semi-algebraic sets}

The study of curvatures of subanalytic sets was started by Fu in [Fu2] . Following ideas of Wintgen [Wi]
and Zaehle [Za], he associated to a compact subanalytic set $X$ a current $NC(X)$, called the normal cycle of $X$, which enables him to define Lipschitz-Killing measures on $X$.

Let $X \subset \mathbb{R}^n$ be a compact subanalytic set. There exists a subanalytic function $f : \mathbb{R}^n \rightarrow
\mathbb{R}$ such that $f \ge 0$, $X =f^{-1}(0)$ and $X_r := f^{-1}([0,r])$ is a smooth compact manifold with boundary for $r$
positive and small. Let $N_r$ be the outside pointing unit normal bundle of $X_r$:
$$N_r = \left\{ (x,v) \ \vert \ x \in f^{-1}(r) \hbox{ and } v =\frac{\nabla f}{\Vert \nabla f \Vert}(x) \right\}.$$
It is a subset of the sphere bundle $S\mathbb{R}^n=\mathbb{R}^n \times S^{n-1}$. Let $D^*(S\mathbb{R}^n)$ be the algebra of
differentials forms on $S\mathbb{R}^n$ and $D_*(S \mathbb{R}^n)$ be the space of all currents with finite mass. The flat
semi-norm on $D_{n-1}(S\mathbb{R}^n)$ is defined by:
$$ \Vert C \Vert^\flat := \left\{ C(\phi) \ \vert \ \phi \in D^{n-1}(S\mathbb{R}^n), \Vert \phi \Vert=1, \Vert d \phi \Vert=1
\right\},$$
where $\Vert - \Vert$ is the comass. The topology on $D_{n-1}(S\mathbb{R}^n)$ induced by $\Vert - \Vert^\flat$ is called the flat
topology. The outside pointing unit normal bundle $N_r$ of $X_r$ defines a current $\tilde{N}_r$ in $D_{n-1}(S\mathbb{R}^n)$.
By the compactness Theorem of Federer-Fleming ([Fe], 4.2.17), there is a sequence of positive real numbers $(r_k)_{k \in \mathbb{N}}$
tending to $0$, such that the sequence of currents $\tilde{N}_{r_k}$ tends to a rectifiable current $\tilde{N}_0$ in
$D_{n-1}(S\mathbb{R}^n)$. Furthermore $\tilde{N}_0$ is a cycle i.e $\partial \tilde{N}_0 =0$ and $\tilde{N}_0 \llcorner \kappa =0$
where $\kappa$ is the canonical $1$-form on $S\mathbb{R}^n$. Then Fu proved a uniqueness theorem which shows that the limit cycle
is independent of the different choices. We call it the normal cycle of $X$ and denote it by $NC(X)$. With this cycle, we can define Lipschitz-Killing measures on Borel sets of $X$.
Let $(x,v) \in S\mathbb{R}^n = \mathbb{R}^n \times S^{n-1}$. Let $\mathcal{B}=(b_1=v,b_2,\ldots,b_n)$ be a direct orthonormal basis of $T_x \mathbb{R}^n$ and let $\mathcal{B}^*=(\sigma_1,\ldots,\sigma_n)$ be its dual basis. Then $(b_2,\ldots,b_n)$ is an orthonormal basis of $T_v S ^{n-1}$, its dual basis being denoted by $\omega_2,\ldots,\omega_n$. For $k \in 
\{0,\ldots,n-1\}$, we set:
$$\psi_k(x,v)= \frac{1}{s_{n-k-1}} \sum_\tau (-1)^{\vert \tau \vert} \omega_{\tau(2)}\wedge \cdots \wedge \omega_{\tau(n-k)} \wedge \sigma_{\tau(n-k-1)} \wedge \cdots 
\wedge \sigma_{\tau(n)},$$ 
where the summation runs over all permutations of $\{2,\ldots,n\}$. Let $\pi : S\mathbb{R}^n \rightarrow S^{n-1}$ be the canonical projection on the unit sphere. For every Borel set $U$ of $X$ and for $k \in \{0,\ldots,n-1\}$, let $\Lambda_k(X,U)$ be defined by:
$$\Lambda_k(X,U)=NC(X) \llcorner \pi^{-1}(U)(\psi_k).$$
For $k=n$, we set $\Lambda_n(X,U) =\mathcal{L}_n(U)$, where $\mathcal{L}_n$ is the $n$-dimensional Lebesgue measure in $\mathbb{R}^n$. 
The measures $\Lambda_k(X,-)$ are called the Lipschitz-Killing measures of $X$. The numbers $\Lambda_k(X,X)$ are the Lipschitz-Killing curvatures of $X$. 
Using the Lipschitz-Killing curvatures, Fu established a Gauss-Bonnet formula and a kinematic formula for compact subanalytic sets (see [Fu2], Theorem 1.8 and Corollary 2.2.2).

In [BK] (see also [BB1] and [BB2]), the authors give a geometric definition of the Lipschitz-Killing measures of closed semi-algebraic sets using stratified Morse theory. Let $X \subset \mathbb{R}^n$ be a closed semi-algebraic set equipped with a semi-algebraic Whitney stratification $\{V_\alpha\}_{\alpha \in \Lambda}$.

 Let us fix a stratum $V$ of $X$ of dimension $d$ with $d<n$. For $x \in V$, we denote by $S_x$ the unit sphere of $(T_x V)^\perp$. For $k \in \{0,\ldots,d\}$, let $\lambda_k^V$ be defined on $V$ by:
$$\lambda_k^V(x)=\frac{1}{s_{n-k-1}} \int_{S_x} \alpha(x,v) \sigma_{d-k} (II_{x,v})dv,$$
where $II_{x,v}$ is the second fundamental form of $V$ in the direction of $v$ and $\sigma_{d-k}$ is its $(d-k)$-th symmetric function. The index $\alpha(x,v)$ is the normal Morse index at $x$ of a function $f : \mathbb{R}^n \rightarrow \mathbb{R}$ such that $f_{\vert X}$ has a Morse critical point at $x$ and $\nabla f(x)=-v$ . This index is defined for almost all $v \in S_x$ (see [GMP, I.1.8] or 
[BK, Lemma 3.5]). For $k \in \{d+1,\ldots,n\}$, we set $\lambda_k^V(x)=0$. 

If $V$ has dimension $n$ then for all $x \in V$, we put $\lambda_0^V(x)=\ldots=\lambda_{n-1}^V(x)=0$ and $\lambda_n^V(x)=1$. 
If $V$ has dimension $0$ then we set:
$$\lambda_0^V(x)=\frac{1}{s_{n-1}} \int_{S^{n-1}} \alpha(x,v) dv,$$
and $\lambda_k^V(x)=0$ if $k>0$.

For every Borel set $U \subset X$ and for every $k \in \{0,\ldots,n\}$, we define $\Lambda'_k(X,U)$ by:
$$\Lambda'_k(X,U) = \sum_{\alpha \in \Lambda} \int_{V_\alpha \cap U} \lambda_k^{V_\alpha} (x) dx.$$
Of course, the measures $\Lambda'_k(X,-)$ coincide with the measures $\Lambda_k(X,-)$ defined by Fu, as explained in [BB2]. From now on, we will keep the unique notation $\Lambda_k(X,-)$ for the Lipschitz-Killing curvatures. 

\section{Gauss-Bonnet formulas for closed semi-algebraic sets}
In this  section, we prove several Gauss-Bonnet type formulas for the Lipschitz-Killing curvatures of a closed semi-algebraic set. 

Let $X \subset \mathbb{R}^n$ be a closed semi-algebraic set. In [Du] Theorem 5.6 and Corollary 5.7, we proved the following Gauss-Bonnet formula:
$$\Lambda_0(X,X)= \chi(X)-\frac{1}{2} \chi(Lk^\infty(X)) 
\hfill -\frac{1}{2 \hbox{vol}(S^{n-1})} 
\int_{S^{n-1}} \chi(Lk^\infty(X \cap v^\perp))dv,$$
where $\Lambda_0(X,X)=\lim_{R \rightarrow +\infty} \Lambda_0(X,X\cap K_R)$ with $(K_R)_{R>0}$ an exhaustive family of compact subsets of $X$ and  $Lk^\infty (X)=X \cap S_R^{n-1}$, $R \gg 1$, is the link at infinity of $X$. 
 
 Here we will express the limits $\lim_{R \rightarrow +\infty} \frac{1}{b_k R^k} \Lambda_k(X, X \cap B_R)$ in terms of the mean-values $\frac{1}{g_n^k} \int_{G_n^k} \chi(Lk^\infty(X\cap H))dH$.

\subsection{The conic case}
Let $X \subset \mathbb{R}^n$ be a closed semi-algebraic cone. In this case for all $H \in G_n^k$, $\chi(Lk^\infty(X \cap H))=\chi(X\cap H \cap S^{n-1})$, and by homogeneity of the Lipschitz-Killing measures: 
$$\frac{1}{R^k} \Lambda_k(X,X \cap B_R)=\Lambda_k(X,X \cap B^n) ,$$
for $k \in \{1,\ldots,n\}$. 
We can equip $X$ with a Whitney stratification: $$X =\sqcup  V_\alpha \sqcup \{0\},$$ where the $V_\alpha$'s are semi-algebraic conic submanifolds. Since $k>0$, we can write:
$$\Lambda_k(X,X \cap B^n)= \sum_{\alpha \in \Lambda} \int_{V_\alpha \cap B^n} \lambda_k^{V_\alpha}(x) dx.$$
Before computing these Lipschitz-Killing integrals, we need some notations and definitions. Let $\tilde{X}$ be defined by $\tilde{X}=X \cap S^{n-1}$. It is a semi-algebraic set of $S^{n-1}$ equipped with the Whitney stratification $\{\tilde{V}\}_{\alpha \in \Lambda} $ where $\tilde{V}_\alpha=V_\alpha \cap S^{n-1}$. Since $S^{n-1}$ is an analytic riemannian manifold of constant sectional curvature equal to $1$, one can associate relative Lipschitz-Killing curvatures with $\tilde{X}$ (see [BB2]). Let $\tilde{V}$ be  a stratum of dimension $\tilde{d}$ with $\tilde{d} < n-1$, let $\tilde{x}$ be  a point in $\tilde{V}$ and let $S_{\tilde{x}}$ be the unit sphere in $(T_{\tilde{x}} \tilde{V}) ^{\perp S^{n-1}}$, the normal space to $\tilde{V}$ at $\tilde{x}$ in $S^{n-1}$. For $k \in \{0,\ldots,\tilde{d}\}$, we set: 
$$\tilde{\lambda}_k^{\tilde{S}} (\tilde{x})=\frac{1}{s_{n-1-k-1}} \int_{S_{\tilde{x}}} \tilde{\alpha}(\tilde{x},\tilde{v}) \sigma_{\tilde{d}-k}(II_{\tilde{x},\tilde{v}}) d\tilde{v},$$
where $II_{\tilde{x},\tilde{v}}$ is the second fundamental form of $\tilde{V}$ in the direction $\tilde{v}$ and $\sigma_{\tilde{d}-k}$ its symmetric function of order $\tilde{d}-k$. The index $\tilde{\alpha}(\tilde{x},\tilde{v})$ is the normal Morse index at $x$ of a function $f :S^{n-1} \rightarrow \mathbb{R}$ such that $f_{\vert \tilde{X}}$ has a Morse critical point at $x$ and $\nabla f(x)=-\tilde{v}$ . For $k \in \{ \tilde{d}+1,\ldots,n-1\}$, we set $\tilde{\lambda}_k^{\tilde{V}}(\tilde{x})=0$. 

If $\tilde{V}$ has dimension $n-1$,  then for all $\tilde{x}$ in $\tilde{V}$, we put  $\tilde{\lambda}_0^{\tilde{V}}(\tilde{x})=\ldots= \tilde{\lambda}_{n-2}^{\tilde{V}}(\tilde{x})=0$ and $\tilde{\lambda}_{n-1}^{\tilde{V}}(\tilde{x})=1$. If $\tilde{V}$ has dimension $0$, then we set:
$$\lambda_0^{\tilde{V}}(x)=\frac{1}{s_{n-2}} \int_{S_{\tilde{x}}} \tilde{\alpha}(\tilde{x},\tilde{v}) dv,$$
and $\lambda_k^{\tilde{V}}(x)=0$ if $k>0$.

For every Borel set $U \subset \tilde{X}$ and for every $k \in \{0,\ldots,n-1\}$, we define:
$$\tilde{\Lambda}_k(\tilde{X},U)= \sum_{\alpha \in \Lambda} \int_{\tilde{V}_\alpha \cap U} \tilde{\lambda}_k^{\tilde{V}_\alpha} (\tilde{x}) d\tilde{x}.$$
Now let us go back to the curvatures $\Lambda_k(X,X\cap B^n)$. Let us a fix a stratum $V$ of dimension $d$ with $d<n$, let $x$ be a point in $V$ and let $\tilde{x}$ in $\tilde{V}$ be defined by $\tilde{x}=\frac{x}{\Vert x \Vert}$. Since $V$ is conic, the normal space at $x$ to $V$ in $\mathbb{R}^n$ is the same as the normal space at $\tilde{x}$ to $\tilde{V}$ in $S^{n-1}$. 
This implies that for almost all $v$ in $(T_x V)^\perp=(T_{\tilde{x}} \tilde{V})^{\perp S^{n-1}}$, $\alpha(x,v)=\tilde{\alpha}(\tilde{x},v)$. Moreover, always by the conic structure, it is not difficult to establish that $\sigma_{d-k} (II_{x,v})=\frac{1}{\Vert x \Vert} \sigma_{d-k} (II_{\tilde{x},v})$ for $k \in \{1,\ldots,d\}$ (see [CM] page 240, for a similar computation). Therefore, we have:
$$\lambda_k^V(x)=\frac{1}{\Vert x \Vert^{d-k}} \tilde{\lambda}_{k-1}^{\tilde{V}} \left( \frac{x}{\Vert x \Vert} \right).$$
If $V$ is of dimension $n$ then for $k\in \{1,\ldots,n-1\}$, $\lambda_k^V(x)=0=\tilde{\lambda}_{k-1}^{\tilde{V}} \left( \frac{x}{\Vert x \Vert} \right)$ and 
$\lambda_n^V(x)=1=\tilde{\lambda}_{n-1}^{\tilde{V}} \left( \frac{x}{\Vert x \Vert} \right)$. 
Integrating $\lambda_k^V$ on $V \cap B^n$ and using spherical coordinates, we obtain:
$$\int_{V \cap B^n} \lambda_k^V(x) dx =\frac{1}{k} \int_{\tilde{V}} \tilde{\lambda}_{k-1}^{\tilde{V}} (\tilde{x}) d\tilde{x},$$ 
and so $\Lambda_k(X,X\cap B^n)=\frac{1}{k} \tilde{\Lambda}_{k-1}(\tilde{X},\tilde{X})$ for $k \in \{1,\ldots,n\}$. 

It remains to relate the relative Lipschitz-Killing curvatures to the mean-values $\frac{1}{g_n^k} \int_{G_n^k} \chi(X \cap H \cap S^{n-1}) dH$. For this, we use the kinematic formula in $S^{n-1}$ (see [Fu1] and [BB2]) and the generalized Gauss-Bonnet formula in the sphere (Theorem 1.2 in [BB2]). For $H \in G_n^k$, $k \in \{1,\ldots,n-1\}$, $X \cap H \cap S^{n-1}= \tilde{X} \cap E_H$ where
$E_H=H  \cap S^{n-1}$ is a $(k-1)$-dimensional sphere of radius $1$. By the generalized Gauss-Bonnet formula, we have:
$$\chi(\tilde{X} \cap E_H)= \sum_{i=0,2,\ldots}^{n-1} \frac{2}{s_i} \tilde{\Lambda}_i(\tilde{X} \cap E_H).$$
Using the notations and the normalizations of [BB2] Theorem 4.1, we can write that:
$$\frac{1}{g_n^k} \int_{G_n^k} \chi(\tilde{X} \cap E_H) dH =\frac{1}{s_{n-1}} \sum_{i=0,2,\ldots}^{n-1} \frac{2}{s_i} \int_G \tilde{\Lambda}_i(\tilde{X},\tilde{X} \cap \sigma E) d\sigma,$$
where $G$ is the isometry group of $S^{n-1}$ and $E$ is a $(k-1)$-dimensional unit sphere in $S^{n-1}$. Since $\tilde{\Lambda}_{k-1}(E,E)=s_{k-1}$ and $\tilde{\Lambda}_i(E,E)=0$ for $i \not= k-1$, we find, applying the kinematic formula, that:
$$\int_G \tilde{\Lambda}_i (\tilde{X},\tilde{X} \cap \sigma E)d\sigma =\frac{s_{n-1}s_i}{s_{n-k+i}}\tilde{\Lambda}_{n-k+i}(\tilde{X},\tilde{X}).$$
Thus for $k \in \{1,\ldots,n\}$, we have:
$$\frac{1}{g_n^k} \int_{G_n^k} \chi(\tilde{X} \cap H) dH =\sum_{i=0,2,\ldots}^{n-1} \frac{2}{s_{n-k+i}} \tilde{\Lambda}_{n-k+i}(\tilde{X},\tilde{X}).$$
When $k=1$, this gives:
$$\frac{1}{g_n^1} \int_{G_n^1} \chi(\tilde{X} \cap H) dH = \frac{2}{s_{n-1}} \tilde{\Lambda}_{n-1}(\tilde{X},\tilde{X}).$$
When $k=2$, this gives:
$$\frac{1}{g_n^2} \int_{G_n^2} \chi(\tilde{X} \cap H) dH = \frac{2}{s_{n-2}} \tilde{\Lambda}_{n-2}(\tilde{X},\tilde{X}).$$
For $k \ge 3$, we have:
$$\frac{1}{g_n^{k-2}} \int_{G_n^{k-2}} \chi(\tilde{X} \cap H) dH =\sum_{i=0,2,\ldots}^{n-1} \frac{2}{s_{n-k+2+i}} \tilde{\Lambda}_{n-k+2+i}(\tilde{X},\tilde{X}),$$
$$= \sum_{i=2,4,\ldots}^{n-1} \frac{2}{s_{n-k+i}} \tilde{\Lambda}_{n-k+i}(\tilde{X},\tilde{X}),$$
and so: 
$$\frac{2}{s_{n-k}} \tilde{\Lambda}_{n-k}(\tilde{X},\tilde{X})=\frac{1}{g_n^k}\int_{G_n^k} \chi(\tilde{X} \cap H)dH -\frac{1}{g_n^{k-2}} \int_{G_n^{k-2}} \chi(\tilde{X} \cap H) dH.$$
Since $\Lambda_{n-k+1}(X,X \cap B^n) =\frac{1}{n-k+1} \tilde{\Lambda}_{n-k}(\tilde{X},\tilde{X})$ and $s_{n-k}=(n-k+1)b_{n-k+1}$, we obtain the following proposition:
\begin{proposition}
Let $X \subset \mathbb{R}^n$ be  a closed semi-algebraic cone and let $\tilde{X}=X \cap S^{n-1}$. Then we have:
$$\frac{\Lambda_n(X,X\cap B^n)}{b_n} = \frac{1}{2g_n^1} \int_{G_n^1} \chi(\tilde{X} \cap H) dH,$$
$$\frac{\Lambda_{n-1}(X,X\cap B^n)}{b_{n-1}} = \frac{1}{2g_n^2} \int_{G_n^2} \chi(\tilde{X} \cap H) dH,$$
and:
$$\displaylines{
\qquad \frac{\Lambda_{n-k}(X,X\cap B^n)}{b_{n-k}} = \frac{1}{2g_n^{k+1}} \int_{G_n^{k+1}} \chi(\tilde{X} \cap H) dH \hfill \cr
\hfill - \frac{1}{2g_n^{k-1}} \int_{G_n^{k-1}} \chi(\tilde{X} \cap H) dH, \qquad}$$
for $k \in \{2,\ldots,n-1\}$, 
or equivalently:
$$\displaylines{
\qquad \frac{\Lambda_{k}(X,X\cap B^n)}{b_{k}} = -\frac{1}{2g_n^{n-k-1}} \int_{G_n^{n-k-1}} \chi(\tilde{X} \cap H) dH \hfill \cr
\hfill + \frac{1}{2g_n^{n-k+1}} \int_{G_n^{n-k+1}} \chi(\tilde{X} \cap H) dH, \qquad}$$
for $k \in \{1,\ldots,n-2\}$. 
\end{proposition}

\subsection{The general case}
We use the same procedure as we did in [Du]. 
We will treat first the case when $X \subset \mathbb{R}^n$ is a closed semi-algebraic of positive codimension. Consider
the set: 
$$X_1 =\left\{ (tx,t) \ \vert \ x \in X, t \in ]0,1] \right\} \subset \mathbb{R}^n \times \mathbb{R}.$$
For each $t \in [0,1]$, put $X_t:= \pi_{\mathbb{R}^n} (\overline{X_1} \cap \pi_\mathbb{R}^{-1}(t))$, where
$\pi_{\mathbb{R}^n} : \mathbb{R}^n \times \mathbb{R} \rightarrow \mathbb{R}^n$ and $\pi_\mathbb{R} :\mathbb{R}^n \times
\mathbb{R} \rightarrow \mathbb{R}$ are the obvious projections. It is easy to see that $X_0$ is a conic semi-algebraic
set. It can be viewed as a tangent cone at infinity of the set $X$. As in Theorem 3.7 of [FM], for each $p \in X_0$,
there is an $\varepsilon_0 >0$ such that if $0<\varepsilon < \varepsilon_0$, then the limit :
$$\phi_0 \cdot 1_{X_t \cap B^n}(p):= \lim_{t \rightarrow 0} \chi(1_{X_t \cap B^n} \cdot 1_{B_\varepsilon(p)})$$
exists and is independent of $\varepsilon$, where $B_\varepsilon(p)$ is the ball of radius $\varepsilon$ centered at $p$.
This resulting function $\phi_0 \cdot 1_{X_t \cap B^n}$ is a semialgebraically constructible function supported on $X_0
\cap B^n$. It is called the specialization of the characteristic function of the family $X_t \cap B^n$ at $t=0$. We will denote it by $\phi_X$.

Since any two sets $X_t$ and $X_{t'}$, where $t$ and $t'$ are different from $0$, are homothetic, it is not difficult to see that the
function $\phi_X$ is conic, i.e $\phi_X(\lambda.x)=\phi_X(x)$ if $x\not=0$ and $\lambda$ is a real different from $0$.
Furthermore, as in Theorem 3.7 [FM], one has :
$$NC(\phi_X)=\lim_{t \rightarrow 0} NC(X_t \cap B^n),$$
where $NC(\phi_X)$ is the normal cycle of the constructible function $\phi_X$ (for the definition of the normal cycle of a
constructible function, see [Be] or [FM]). 

For any semi-algebraically constructible function $\psi$, let $Lk^\infty(\psi)$ be the function $1_{S^{n-1}} \cdot \psi$.
We know that $Lk^\infty(\phi_X)=\phi_0 \cdot 1_{X_t \cap S^{n-1}}$, the specialization at $0$ of the characteristic function of the family $X_t \cap S^{n-1}$ (see [Du], Lemma 5.3). 
We need auxiliary lemmas.
\begin{lemma}
Let $S \subset \mathbb{R}^n$ be a semi-algebraic smooth manifold. There exists a semi-algebraic set $\Sigma_k' \subset G_n^k$,
 $k \in \{1,\ldots,n-1\}$, of positive codimension such that if $H \notin \Sigma_k'$, $H$ intersects  $S \setminus \{ 0\}$ transversally.
\end{lemma}

{\it Proof.} Assume that $S$ has dimension $d$ and fix $k \in \{1,\ldots,n-1\}$. Let $W$ be defined by:
$$\displaylines{
\qquad W= \big\{ (x,v_1,\ldots,v_{n-k}) \in \mathbb{R}^n \times (\mathbb{R}^n)^{n-k} \ \vert  \hfill \cr
\hfill   x \in S \setminus \{ 0\} \hbox{ and } 
\langle x,v_1 \rangle = \cdots = \langle x, v_{n-k} \rangle=0 \big\}. \qquad}
$$
It is a smooth semi-algebraic manifold of dimension $d+(n-1)(n-k)$. Let $\pi_2$ be the following projection:
$$\pi_2 : W \rightarrow (\mathbb{R}^n)^{n-k}, (x,v_1,\ldots,v_{n-k}) \mapsto (v_1,\ldots,v_{n-k}).$$
 Bertini-Sard's theorem implies that the set of critical values of $\pi_2$ is a semi-algebraic set of positive codimension. $\hfill \Box$

Now we can compare the sets $(X \cap H)_0$ and $X_0 \cap H$. First, let us equip $\overline{X_1}$ with a Whitney stratification compatible with $X_0 \times \{0\}$. This induces a Whitney stratification $\Sigma$ on $X_0$. 
\begin{lemma}
Let $k \in \{1,\ldots,n-1\}$ and let $H \in G_n^k$. If $H$ intersects $X_0$ equipped with the stratification $\Sigma$ transversally then $(X \cap H)_0=X_0 \cap H$.
\end{lemma}

{\it Proof.} The inclusion $(X\cap H)_0 \subset X_0 \cap H$ is clear because $\overline{X_1 \cap H} \subset \overline{X_1} \cap H \times \mathbb{R}$. 
If $X_0 \cap H = \{0\}$ then the result is trivial. If $X_0 \cap H \not= \{0\}$ then let us choose a point $x \not= 0$ in $X_0 \cap H$. If $S$ is the stratum that contains $x$ then $H \times \mathbb{R}$ intersects $S\times \{0\}$ transversally in $\mathbb{R}^n \times \mathbb{R}$. By the Whitney condition $(a)$, it intersects also transversally around $x$ the strata that contains $S \times \{0\}$ in their closure. Since $(x,0)$ belongs to $\overline{X_1}$, there exists at least one stratum $T \subset X_1$ that contains $S \times \{0\}$ in its closure. Hence we can find a sequence of points in $T \cap H$ of the form $(t_n y_n,t_n)$, $y_n \in X$, $t_n \in ]0,1]$, that converges to $(x,0)$. Since $H$ is a linear subspace, $y_n \in H$ and so $x \in (X \cap H)_0$. $\hfill \Box$

For any $k \in \{1,\ldots,n-1\}$, for any $H \in G_n^k$ and for any semi-algebraically constructible function $\psi$, we denote by $\psi_H$ the semi-algebraically constructi\-ble function $\psi \cdot 1_H$.
\begin{lemma}
Let $k \in \{1,\ldots,n-1\}$ and let $H \in G_n^k$. If $H$ intersects $X_0$ equipped with $\Sigma$  transversally then $Lk^\infty \left[ (\phi_X)_H \right]$ is the specialization of the family $H \cap X_t \cap S^{n-1}$ at $0$.
\end{lemma}

{\it Proof.} We prove the lemma by induction. For $k=n-1$, this lemma is true by Lemmas 5.3 and 5.4 in [Du]. Let us assume that it is true for $k \in \{2,\ldots,n-1\}$ and let $H \in G_n^{k-1}$ be a $(k-1)$-plane that intersects $X_0$ transversally. We have to prove that:
$$\lim_{t \rightarrow 0} \chi (X_t \cap B^n \cap B_\varepsilon(p)) =\lim_{t \rightarrow 0} \chi (X_t \cap H \cap S^{n-1} \cap B_\varepsilon(p)) ,$$
for every $p$ in $X_0 \cap H \cap S^{n-1}$. 
If $H \cap X_0 \setminus \{0\}$ is empty then the result is true. If $H \cap X_0 \setminus \{0\} \not= \emptyset$ then let $L$ be a $k$-plane that contains $H$. Then $L$ intersects $X_0 $ transversally and so, by the induction hypothesis, we have:
$$\lim_{t \rightarrow 0} \chi (X_t \cap B^n \cap B_\varepsilon(p)) =\lim_{t \rightarrow 0} \chi (X_t \cap L \cap S^{n-1} \cap B_\varepsilon(p)) ,$$
for every $p$ in $X_0 \cap L \cap S^{n-1}$.  Let $v$ be an unit vector in $L$ such that $H=L \cap v^\perp$. Applying Lemma 5.4 in [Du] to $X \cap L$, we find that:
$$\lim_{t \rightarrow 0} \chi ((X \cap L)_t \cap B^n \cap B_\varepsilon(p)) =\lim_{t \rightarrow 0} \chi ((X\cap L)_t \cap v^\perp \cap S^{n-1} \cap B_\varepsilon(p)) ,$$
for every $p$ in $(X\cap L)_0 \cap v^\perp \cap S^{n-1}$. We conclude with the fact that for $t \not= 0$, $(X \cap L)_t=X_t \cap L$ and the fact that $(X \cap L)_0=X_0 \cap L$, by the previous lemma. $\hfill \Box$

Combining this lemma with Lemma 3.2, we get:
\begin{corollary}
Let $k \in \{1,\ldots,n\}$. For almost all $H \in G_n^k$, $Lk^\infty \left[ (\phi_X)_H \right]$ is the specialization of the family $H \cap X_t \cap S^{n-1}$ at $0$.
\end{corollary}

\begin{proposition}
We have, for $k \in \{1,\ldots,n-2\}$:
$$\displaylines{
\qquad \frac{1}{b_k} NC(\phi_X) \llcorner \pi^{-1}(\mathring{B^n}) (\psi_k) =
 -\frac{1}{2g_n^{n-k-1}} \int_{G_n^{n-k-1}} \chi(LK^\infty[(\phi_X)_H]) dH \hfill \cr
\hfill + \frac{1}{2g_n^{n-k+1}} \int_{G_n^{n-k+1}} \chi(LK^\infty[(\phi_X)_L]) dL, \qquad \cr
}$$
and:
$$\frac{1}{b_{n-1}} NC(\phi_X) \llcorner \pi^{-1} (\mathring{B^n}) (\psi_{n-1}) =
  \frac{1}{2g_n^{2}} \int_{G_n^{2}} \chi(LK^\infty[(\phi_X)_H]) dH,$$
$$\frac{1}{b_{n}} NC(\phi_X) \llcorner \pi^{-1} (\mathring{B^n}) (\psi_{n}) =
  \frac{1}{2g_n^{1}} \int_{G_n^{1}} \chi(LK^\infty[(\phi_X)_H]) dH,$$
  where $\mathring{B^n}$ is the interior of $B^n$.
 \end{proposition}
 {\it Proof.} We prove the equality for $k \in \{1,\ldots,n-2\}$. Since the two sides of the equality are linear over $\mathbb{Z}$ in the constructible function $\phi_X$ supported on $B^n$, it is enough to establish the equality for the function $1_{X \cap B^n}$ where $X$ is a closed semi-algebraic cone. In this case, we have to prove:
 $$\displaylines{
 \qquad \frac{1}{b_k} NC(X \cap B^n) \llcorner \pi^{-1}(\mathring{B^n}) (\psi_k) =
 -\frac{1}{2g_n^{n-k-1}} \int_{G_n^{n-k-1}} \chi(LK^\infty(X \cap H)) dH \hfill \cr
\hfill + \frac{1}{2g_n^{n-k+1}} \int_{G_n^{n-k+1}} \chi(LK^\infty(X \cap L)) dL. \qquad \cr
}$$
This is exactly the formula proved in Proposition 3.1. $\hfill \Box$

\begin{theorem}
Let $X \subset \mathbb{R}^n$ be a closed semi-algebraic set of positive codimension. For $k \in \{1,\ldots,n-2\}$, we have:
$$\displaylines{
\qquad \lim_{R \rightarrow + \infty} \frac{\Lambda_k(X,X\cap B_R)}{b_k R^k} =- \frac{1}{2g_n^{n-k-1}} \int_{G_n^{n-k-1}} \chi(Lk^\infty(X \cap H)) dH \hfill \cr
\hfill + \frac{1}{2g_n^{n-k+1}} \int_{G_n^{n-k+1}} \chi(Lk^\infty(X \cap L)) dL. \qquad \cr
}$$
and:
$$\lim_{R \rightarrow + \infty} \frac{\Lambda_{n-1}(X,X\cap B_R)}{b_{n-1} R^{n-1}} = \frac{1}{2g_n^2} \int_{G_n^2} \chi(Lk^\infty(X \cap H)) dH,$$
$$\lim_{R \rightarrow + \infty} \frac{\Lambda_{n}(X,X\cap B_R)}{b_{n} R^{n}} = \frac{1}{2g_n^1} \int_{G_n^1} \chi(Lk^\infty(X \cap H)) dH.$$
\end{theorem}
{\it Proof.} The last equality is clearly true because by the assumption on the codimension of $X$, both sides vanish.

Let $k \in \{1,\ldots,n-1\}$. First remark that $\Lambda_k(X,X\cap B_R)=\Lambda_k(X,X\cap \mathring{B_R})$.  On the one hand, the homogeneity of the Lipschitz-Killing measures implies that:
$$\frac{\Lambda_k(X,X\cap \mathring{B_R})}{ R^k} = \Lambda_k(X_{\frac{1}{R}}, X_{\frac{1}{R}} \cap \mathring{B^n}) =
NC(X_{\frac{1}{R}} \cap B^n) \llcorner \pi^{-1}(\mathring{B^n})(\psi_k).$$
Since $NC (\phi_X)=\lim_{R \rightarrow + \infty} NC(X_{\frac{1}{R}} \cap B^n)$, we obtain that:
$$\lim_{R \rightarrow + \infty} \frac{\Lambda_k(X,X\cap B_R)}{b_k R^k} =\frac{1}{b_k} NC(\phi_X) \llcorner \pi^{-1}(\mathring{B^n})(\psi_k).$$
On the other hand, for almost all $H$ in $G_n^{n-k-1}$, $Lk^\infty(\phi_X)_H$ is the specialization at $0$ of the characteristic function of the family $X_t \cap H \cap S^{n-1}$ at $0$. Therefore: 
$$\chi(Lk^\infty(\phi_X)_H)=NC(Lk^\infty(\phi_X)_H)(\psi_0)=\lim_{t \rightarrow 0} NC(X_t \cap H \cap S^{n-1})(\psi_0)=$$ 
$$\lim_{t \rightarrow 0} \chi(X_t \cap H \cap S^{n-1})=\chi(Lk^\infty(X \cap H)).$$ $\hfill \Box$

\begin{corollary}
Let $X \subset \mathbb{R}^n$ be a closed semi-algebraic set. For $k \in \{1,\ldots,n-2\}$, we have:
$$\displaylines{
\qquad \lim_{R \rightarrow + \infty} \frac{\Lambda_k(X,X\cap B_R)}{b_k R^k} =- \frac{1}{2g_n^{n-k-1}} \int_{G_n^{n-k-1}} \chi(Lk^\infty(X \cap H)) dH \hfill \cr
\hfill + \frac{1}{2g_n^{n-k+1}} \int_{G_n^{n-k+1}} \chi(Lk^\infty(X \cap L)) dL. \qquad  \cr
}$$
and:
$$\lim_{R \rightarrow + \infty} \frac{\Lambda_{n-1}(X,X\cap B_R)}{b_{n-1} R^{n-1}} = \frac{1}{2g_n^2} \int_{G_n^2} \chi(Lk^\infty(X \cap H)) dH,$$
$$\lim_{R \rightarrow + \infty} \frac{\Lambda_{n}(X,X\cap B_R)}{b_{n} R^{n}} = \frac{1}{2g_n^1} \int_{G_n^1} \chi(Lk^\infty(X \cap H)) dH.$$
\end{corollary}
{\it Proof.} It remains to prove the case when $X$ has dimension $n$. We proceed as in [Du], Corollary 5.7. Let $i : \mathbb{R}^n \rightarrow \mathbb{R}^{n+1}$ be the natural embedding of 
$\mathbb{R}^n$ in $\mathbb{R}^{n+1}$. For $k \in \{1,\ldots,n\}$, let us denote by $\Lambda_k^n$ (respectively $\Lambda_k^{n+1}$) the Lipschitz-Killing measure of $X$ as a semi-algebraic set in $\mathbb{R}^n$ (respectively $\mathbb{R}^{n+1}$). Since $i$ is a semi-algebraic isometry, by Theorem 5.0 in [Fu2] or Proposition 9.2 in [BK], one has:
$$\Lambda_k^n(X,X \cap B_R^n)=\Lambda_k^{n+1}(X,X \cap B_R^{n+1}).$$
Let us assume that $k \in \{1,\ldots,n-2\}$ and let us apply the previous theorem to $\Lambda_k^{n+1}$. We find:
$$\displaylines{
\qquad \lim_{R \rightarrow + \infty} \frac{\Lambda_k^{n+1}(X,X\cap B_R)}{b_k R^k} =- \frac{1}{2g_n^{n-k}} \int_{G_n^{n-k}} \chi(Lk_{n+1}^\infty(X \cap H)) dH \hfill \cr
\hfill + \frac{1}{2g_n^{n-k+2}} \int_{G_n^{n-k+2}} \chi(Lk_{n+1}^\infty(X \cap L)) dL, \qquad \cr
}$$
where $Lk_{n+1}^\infty(-)$ denotes the link at infinity in $\mathbb{R}^{n+1}$. It is clear that:  $$Lk^\infty_{n+1}(X)=Lk_n^\infty(X)=Lk^\infty(X),$$ because $X \subset \mathbb{R}^n$. Let us compare now $\frac{1}{g_n^q} \int_{G_n^q} \chi(Lk^\infty(X \cap H))dH$ and $\frac{1}{g_{n+1}^{q+1}} \int_{G_{n+1}^{q+1} } \chi(Lk^\infty(X \cap H))dH$ for $0<q<n$. We can write:
$$\int_{G_n^q} \chi(Lk^\infty(X \cap H))dH= \sum_{i=1}^m \chi_i \hbox{vol}(A_i),$$
where the $A_i$'s are open semi-algebraic sets of $G_n^q$. Therefore:
$$\int_{G_{n+1}^{q+1}} \chi(Lk^\infty(X \cap H))dH= \sum_{i=1}^m \chi_i \hbox{vol}(\tilde{A}_i),$$
where $\tilde{A}_i = \{ L \in G_{n+1}^{q+1} \ \vert \ L \cap \mathbb{R}^n \in A_i \}$. We conclude remarking that vol$(A_i)=\frac{\hbox{vol}(G_n^q)}{\hbox{vol}(G_{n+1}^{q+1})} \cdot \hbox{vol}(\tilde{A}_i)$. $\hfill \Box$

Combining Corollary 3.9 with Corollary 5.7 of [Du], we obtain the following Gauss-Bonnet type formula which relates the Euler characteristic of a closed semi-algebraic set to the Lipschitz-Killing curvatures.
\begin{theorem}
Let $X \subset \mathbb{R}^n$ be a closed semi-algebraic set. We have:
$$\chi(X)= \Lambda_0(X,X) + \sum_{k=1}^n    \lim_{R \rightarrow + \infty} \frac{\Lambda_k(X,X \cap B_R)}{b_k R^k}.$$
\end{theorem}
$\hfill \Box$

It is worth remarking that in all these results, the choice of the origin in $\mathbb{R}^n$ does not matter. Namely if $x_0$ is a point in $\mathbb{R}^n$, $B_R(x_0)$ (resp. $S_R(x_0)$) is the ball $x_0 +B_R$ (resp. the sphere $x_0+S_R)$,
$H_{x_0}$ is the affine space $x_0+H$ and $Lk^\infty_{x_0}(X)= X \cap S_R(x_0)$ with $R \gg1$, then the following formulas hold:
$$\displaylines{
\qquad  \Lambda_0(X,X)= \chi(X)-\frac{1}{2} \chi(Lk_{x_0}^\infty(X)) \hfill \cr
\hfill -\frac{1}{2g_n^{n-1}} 
\int_{G_n^{n-1}} \chi(Lk_{x_0}^\infty(X \cap H_{x_0} ))dH, \qquad \cr
}$$
and for $k \in \{1,\ldots,n-2\}$:
$$\displaylines{
\qquad \lim_{R \rightarrow + \infty} \frac{\Lambda_k(X,X\cap B_R(x_0))}{b_k R^k} =- \frac{1}{2g_n^{n-k-1}} \int_{G_n^{n-k-1}} \chi(Lk_{x_0}^\infty(X \cap H_{x_0})) dH \hfill \cr
\hfill + \frac{1}{2g_n^{n-k+1}} \int_{G_n^{n-k+1}} \chi(Lk_{x_0}^\infty(X \cap L_{x_0})) dL. \qquad  \cr
}$$
and:
$$\lim_{R \rightarrow + \infty} \frac{\Lambda_{n-1}(X,X\cap B_R(x_0))}{b_{n-1} R^{n-1}} = \frac{1}{2g_n^2} \int_{G_n^2} \chi(Lk_{x_0}^\infty(X \cap H_{x_0})) dH,$$
$$\lim_{R \rightarrow + \infty} \frac{\Lambda_{n}(X,X\cap B_R(x_0))}{b_{n} R^{n}} = \frac{1}{2g_n^1} \int_{G_n^1} \chi(Lk_{x_0}^\infty(X \cap H_{x_0})) dH.$$
Therefore, we also have:
$$\chi(X)= \Lambda_0(X,X) + \sum_{k=1}^n    \lim_{R \rightarrow + \infty} \frac{\Lambda_k(X,X \cap B_R(x_0))}{b_k R^k}.$$

\section{Application to the smooth case}

Let us recall first some results about complete manifolds. The first one is the well-known Cohn-Vossen inequality [Co]. It states that if $(M,g)$ is a complete connected oriented riemannian surface of finite topological type and with absolute integral Gauss curvature $K$ then :
$$2 \pi \chi(M) -\int_M Kdx \ge 0.$$
In [Sh], Shiohama made this phenomena more precise and gave a geometric expression for the Gauss-Bonnet defect. Namely he proved that:
$$2 \pi \chi(M) -\int_M Kdx = \lim_{t \rightarrow + \infty} \frac{L(t)}{t} =  \lim_{t \rightarrow + \infty} \frac{A(t)}{t^2},$$
where $L(t)$ denotes the length of the geodesic distance circle in distance $t$ and $A(t)$ denotes the area of the geodesic disc with radius $t$, the center of the disc being arbitrary.

For higher dimensional complete manifolds, the situation is much less understood. In [DK], Dillen and Kuehnel study submanifolds of $\mathbb{R}^n$ that they call submanifolds with cone-like ends (see Definition 5.3 in [DK]). Roughly speaking, a submanifold $M$ in $\mathbb{R}^n$ is with cone-like ends if it admits a finite number of ends and each end is ``close" to an open subset of a cone. One can say that somehow it behaves at infinity like a semi-algebraic set. Moreover one can associate with such a manifold $M$ a submanifold $M_\infty \subset S^{n-1}$ called the set of limit directions. It can be viewed as a kind of link at infinity of the manifold $M$. If $d$ is the dimension of $M$, then $M_\infty$ has dimension $d-1$. Dillen and Kuehnel proved that if $M$ is complete then: 
$$ \chi(M) -\frac{1}{s_{n-1}} \int_M K_0 dx =
\sum_{0 \le 2i \le d-1} \frac{1}{s_{n-d+2i-1}s_{d-1-2i} } \int_{M_\infty} \tilde{K}_{2i}d \tilde{x},$$ 
where $K_0$ is the Lipschitz-Killing curvature of $M$ and the $\tilde{K}_{2i}$ are the spherical Lipschitz-Killing curvature of $M_\infty$ (see [So] for an hyperbolic version of this result). We refer to [DK] for an account on Gauss-Bonnet formulas for complete manifolds. 

In the sequel, we will apply the results of Section 3 to the smooth case to get new Gauss-Bonnet formulas for smooth closed semi-algebraic sets.
Let $X \subset \mathbb{R}^n$ be a closed semi-algebraic set which is a manifold of dimension $d$, $1 \le d \le n-1$. In this situation, there is only one stratum $X$ and for each $k \in \{0,\ldots,d\}$ and each $x \in X$:
$$\lambda_k^X(x)=\frac{1}{s_{n-k-1}} \int_{S_x} \sigma_{d-k}(II_{x,v}) dv,$$
where $II_{x,v}$ is the second fundamental form of $X$ in the direction of $v$ and $\sigma_{d-k}$ is its $(d-k)$-th symmetric function. Furthermore for $k>d$, $\lambda_k^X(x)=0$. 

Let us denote by $K_i(x)$ the real number $\int_{S_x} \sigma_i(II_{x,v})dv$, $i \in \{0,\ldots,d\}$. We remark that $K_i(x)=0$ if $i$ is odd because $\sigma_i(II_{x,-v})=(-1)^i \sigma_i(II_{x,v})$. 
Moreover $K_d(x)=2LK(x)$ where $LK$ is defined in [Du] or [La]. 
The $K_i$'s are the classical Lipschitz-Killing-Weyl curvatures. When $i$ is even, $\frac{1}{s_{n-d+i-1}} K_i$ is an intrinsic quantity, namely it depends only on the inner metric on $X$ (see [DK]). With these notations, Theorem 3.7 can be restated in the following form:
\begin{theorem}
Let $X \subset \mathbb{R}^n$ be a closed semi-algebraic set which is a smooth submanifold of dimension $d$, $1 \le d \le n-1$. We have, for $i \in \{1,\ldots,d-1\}$:
$$\displaylines{
\lim_{R \rightarrow + \infty} \frac{\int_{X \cap B_R} K_i dx}{s_{n-d+i-1} b_{d-i} R^{d-i}} =- \frac{1}{2g_n^{n-d+i-1}} \int_{G_n^{n-d+i-1}} \chi(Lk^\infty(X \cap H)) dH \cr
\hfill + \frac{1}{2g_n^{n-d+i+1}} \int_{G_n^{n-d+i+1}} \chi(Lk^\infty(X \cap L)) dL, \qquad
}$$
and:
$$\lim_{R \rightarrow + \infty} \frac{\hbox{\em vol}(X\cap B_R)}{b_{d} R^{d}} = \frac{1}{2g_n^{n-d+1}} \int_{G_n^{n-d+1}} \chi(Lk^\infty(X \cap H)) dH.$$
\end{theorem}
{\it Proof.} We apply Theorem 3.7. For $i=0$, we remark that for almost all $H \in G_n^{n-d-1}$, the set $X \cap H$ is empty. $\hfill \Box$

Note that when $i$ is odd, the equality is trivial because for almost all $L \in G_n^{n-d+i-1}$ (resp. $H \in G_n^{n-d+i+1}$), $Lk^\infty(X \cap L)$ (resp. $Lk^\infty(X \cap H)$) is a compact odd dimensional manifold (or empty). Since $X$ is smooth, $X \cap H$ is generically smooth and when it is odd dimensional, $\chi(Lk^\infty(X \cap H))=2 \chi(X \cap H)$. Hence, we can reformulate Theorem 4.1 as follows: 
\begin{theorem}
Let $X \subset \mathbb{R}^n$ be a closed semi-algebraic set which is a smooth submanifold of dimension $d$, $1 \le d \le n-1$. We have, for $i \in \{1,\ldots,d-1\}$ and $i$ even:
$$\displaylines{
\qquad \lim_{R \rightarrow + \infty} \frac{\int_{X \cap B_R} K_i dx}{s_{n-d+i-1} b_{d-i} R^{d-i}} =- \frac{1}{g_n^{n-d+i-1}} \int_{G_n^{n-d+i-1}} \chi(X \cap H) dH  \hfill \cr
\hfill + \frac{1}{g_n^{n-d+i+1}} \int_{G_n^{n-d+i+1}} \chi(X \cap L) dL, \qquad \cr
}$$
and:
$$\lim_{R \rightarrow + \infty} \frac{\hbox{\em vol}(X\cap B_R)}{b_{d} R^{d}} = \frac{1}{g_n^{n-d+1}} \int_{G_n^{n-d+1}} \chi(X \cap H) dH.$$
\end{theorem}
$\hfill \Box$

It is interesting to see that if we apply this theorem to the case $d$ odd and $i=d-1$, we get the following equality:
$$\chi(X)=\lim_{R \rightarrow + \infty} \frac{\int_{X \cap B_R} K_{d-1}dx}{s_{n-2} b_{1} R} + \frac{1}{g_n^{n-2}} \int_{G_n^{n-2}} \chi(X \cap H) dH,$$
which is a kind of Gauss-Bonnet formula for odd-dimensional closed semi-algebraic manifolds. 

When we apply Theorem 3.9, we get another Gauss-Bonnet type formula for closed semi-algebraic submanifolds.
\begin{theorem}
Let $X \subset \mathbb{R}^n$ be a closed semi-algebraic set which is a smooth submanifold of dimension $d$, $1 \le d \le n-1$. If $d$ is even, we have:
$$\chi(X)= \frac{1}{s_{n-1}} \int_X K_d(x)dx + \sum_{i=0}^{\frac{d-2}{2}} \lim_{R \rightarrow + \infty} \frac{1}{s_{n-d+2i-1}b_{d-2i}R^{d-2i}} \int_{X \cap B_R} K_{2i}dx.$$
If $d$ is odd, we have:
$$\chi(X)=  \sum_{i=0}^{\frac{d-1}{2}} \lim_{R \rightarrow + \infty} \frac{1}{s_{n-d+2i-1}b_{d-2i}R^{d-2i}} \int_{X \cap B_R} K_{2i}dx.$$
\end{theorem}
$\hfill \Box$

As in the general case, the choice of the origin does not matter.
This theorem gives an answer to Question 2, p.197 in [DK] for closed semi-algebraic submanifolds of the euclidian space. It is also a generalization of Shiohama's formula. 
It would be interesting to enlarge the class of submanifolds for which this Gauss-Bonnet formula is satisfied.

\end{document}